\documentclass[12pt,a4paper]{amsart}
\usepackage{amssymb,latexsym,euscript,amscd,amsmath}

\input xypic

\DeclareMathAlphabet{\E}{U}{eus}{m}{n}
\newcommand{\Spec}{\mathrm{Spec}}

\newcommand{\Hilb}{\ensuremath{\mathrm{Hilb}}}
\newcommand{\p}{\ensuremath{\mathfrak{p}}}
\renewcommand{\H}{\ensuremath{\mathrm{H}}}

\renewcommand{\O}{{\E O}}

\theoremstyle{plain}

\newtheorem{thm}{Theorem}
\newtheorem{lem}{Lemma}
\newtheorem{prop}{Proposition}
\newtheorem{cor}{Corollary}

\theoremstyle{remark}
\newtheorem{remark}{Remark}

\theoremstyle{definition}

\title[Frobenius splitting of Hilbert schemes]
{Frobenius splitting of Hilbert schemes of points on
surfaces}

\author{Shrawan Kumar, Jesper Funch Thomsen }

\address
{Department of Mathematics, University of North Carolina, Chapel Hill,
NC 27599-3250, USA and Matematisk Institut, Aarhus Universitet,
Ny Munkegade, DK-8000 \AA rhus C, Denmark.}
\email{kumar@math.unc.edu, funch$@$imf.au.dk}

\begin{document}

\bigskip

\begin{abstract}
Let $X$ be a quasiprojective smooth surface defined over an
algebraically
closed field of positive characteristic. In this note we show
that if $X$ is Frobenius split then so is the Hilbert scheme $\Hilb^n(X)$
of $n$ points in $X$. In particular, we get the higher
cohomology vanishing for ample line bundles on $\Hilb^n(X)$
when $X$ is projective and Frobenius split.
\end{abstract}

\bigskip

\maketitle
{\bf Introduction}
\vskip2ex

Let $X$ be a quasiprojective smooth surface defined over an
algebraically
closed field $k$ of positive characteristic $p$. For an integer $n \geq 1$,
let $X^{(n)}$ be the
$n$-th  symmetric product of $X$ and let $X^{[n]}$ denote the 
 Hilbert
scheme  of $n$ points in $X$ (parametrizing the
zero dimensional closed subschemes of $X$ of length $n$).  Recall that 
$X^{[n]}$ is smooth and there 
is a birational 
 `Hilbert-Chow' morphism
 $\psi : X^{[n]} \rightarrow X^{(n)}$,
which to each zero dimensional closed subscheme in $X$ of length $n$ associates
its support (with multiplicities).  
Let $X^{(n)}_{*}$ denote the open locus of $X^{(n)}$ corresponding to
the set of n-tuples with at least $n-1$ distinct points and let 
 $X^{[n]}_{*}$ denote its inverse image under $\psi$. We show that 
$\psi : X^{[n]}_{*} \to X^{(n)}_{*}$ is a crepant resolution if $p > 2$,
 in the sense 
that  $X^{(n)}_*$ is Gorenstein
such that its dualizing line bundle $\omega_{X^{(n)}_*}$ 
 pulls back to the canonical bundle $\omega_{X^{[n]}_*}$
on $X^{[n]}_*$ under $\psi$ (cf. Theorem \ref{crepant}). In fact, if $ 
p > n $, 
$\psi : X^{[n]} \to X^{(n)}$ itself is a crepant resolution (cf. Corollary 
\ref{beauville}).
(This generalizes the corresponding result in char. $0$ due to Beauville.)
We make crucial use of our Theorem 1 to prove the following main result of this
paper: 

Let $X$ be as above and  $p > 2$ .
Then, for any $n \geq 1,$  the  Hilbert scheme $X^{[n]} $ is Frobenius split
(cf. Theorem \ref{main}).
In particular, if $X$, in addition, is projective and $L$ is an ample line 
bundle on  $X^{[n]} $, then $L$ has vanishing higher cohomology (cf. 
Corollary \ref{amplevanishing}).  

The contents of the paper are as follows: Section 1 is devoted to recalling 
the definition of Hilbert schemes, and Section 2 is devoted  to the basic
definitions of Frobenius splitting. Sections 3 and 4 are devoted 
to  proving that $\psi$ is a crepant resolution. We prove our main theorem 
(Theorem \ref{main}) in  Section 5.

We thank M. Brion, S. Hansen and V. Mehta for some helpful conversations. 
Part of this 
work was done while the first author was visiting  Ecole Normale 
Sup\'erieure, Paris,  hospitality of which is gratefully acknowledged. 
The second author is pleased to 
thank UNC for its hospitality during his visit.
The second author was partially 
supported by the TMR programme ``Algebraic Lie Representations'' 
(ECM Network Contract No.~ERB FMRX-CT 97/0100).

\maketitle

\section{Hilbert schemes of points}
Let $X$ be a quasiprojective variety defined over an algebraically
closed field $k$. Fix an integer $n \geq 1 $. The Hilbert
scheme $X^{[n]}= \Hilb^n(X)$ of $n$ points in $X$ parameterizes
zero dimensional closed subschemes of $X$ of length $n$.
The scheme $\Hilb^n(X)$ is quasiprojective and in fact projective
when $X$ is so.

\subsection{Symmetric products}
\label{symmetric}
Let $X^n=X \times \dots \times X$ denote the n-fold product
of $X$, and let $S_n$ denote the symmetric group on $n$
letters. Then $S_n$ acts on $X^n$ by permuting the factors.
As $X^n$ is quasiprojective and $S_n$ is finite, the
geometric quotient of this action exists (cf. \cite{Serre},
 Chap. III, \S 14).
The quotient is denoted by $X^{(n)}$ and is called the
$n$-th {\it symmetric product} of $X$. Let $\Phi : X^n \rightarrow
X^{(n)}$
denote the quotient map.

Points in $X^{(n)}$ correspond to unordered tuples
of (not necessarily distinct) $n$ points in $X$. The open subset of
$X^{(n)}$
consisting of the
 tuples of $n$ distinct points is denoted by $X^{(n)}_{**}$. If $X$ is
smooth,
the variety $X^{(n)}$ is smooth along $X^{(n)}_{**}$ and moreover it is
singular along the
complement of $X^{(n)}_{**}$ if dim $X \geq 2$ (cf. \cite{Fogarty2}, \S 2). Clearly,  the codimension
of $X^{(n)} \setminus X^{(n)}_{**}$ in $X^{(n)}$ is equal to dim $X$.
Let $X^{(n)}_{*}$ denote the open locus of $X^{(n)}$ corresponding to
the set of n-tuples with at least $n-1$ distinct points.

\subsection{Hilbert-Chow morphism (\cite{Fogarty1}, \S 2)}

Let $X^{[n]}_{red}$ denote the underlying reduced
subscheme of $X^{[n]}$. The {\it Hilbert-Chow morphism}
is the map $\psi : X^{[n]}_{red} \rightarrow X^{(n)}$,
which to each zero dimensional closed subscheme in $X$ of length $n$ associates
its support (with multiplicities). The Hilbert-Chow
morphism is birational, being an isomorphism over
the open set $X^{(n)}_{**}$.

When $X$ is a smooth surface, the Hilbert scheme
$X^{[n]}$ is also smooth (in particular reduced). Hence,
in this case, $\psi$
is a desingularization of the symmetric product
$X^{(n)}$.

\section{Frobenius splitting - basic definitions }
Let
$\pi:X\rightarrow \Spec(k)$ be a scheme defined over an algebraically
closed field $k$ of positive characteristic $p$. The
{\it absolute Frobenius morphism} on $X$ is the identity on point
spaces and raising to the $p$-th power locally on functions. The
absolute
Frobenius morphism is {\it not} a morphism of $k$-schemes.
Let $X'$ be the scheme obtained from $X$ by base change with the
absolute Frobenius morphism on $\Spec(k)$, i.e., the underlying
topological
space of $X'$ is that of $X$ with the same structure sheaf $\O_X$ of
rings,
only
the underlying $k$-algebra structure on $\O_{X'}$ is twisted as
$\lambda\odot f = \lambda^{1/p} f$, for $\lambda\in k$ and $f\in
\O_{X'}$.
Using this description of $X'$,
{\it the relative Frobenius
morphism} $F: X\rightarrow X'$ is defined in the same way as the
absolute Frobenius morphism and it is a morphism of $k$-schemes.

\subsection{Frobenius splitting \cite{MehtaRamanathan}}

Recall that  a variety $X$ is called
{\it Frobenius split} if the homomorphism $\O_{X'}\rightarrow F_*\O_X$
of
$\O_{X'}$-modules is split.
A homomorphism $\sigma:F_*\O_X\rightarrow \O_{X'}$ is a splitting
of $\O_{X'}\rightarrow F_*\O_X$
(called a {\it Frobenius splitting}) if and only if $\sigma(1)=1$.

When $X$ is a smooth variety with canonical bundle $\omega_X$,
there is a natural isomorphism of $\O_{X'}$-modules
(\cite{MehtaRamanathan}):
$$F_* (\omega_X^{1-p}) \cong Hom_{\O_{X'}}(F_* \O_X, \O_{X'}).$$
In this way global sections of $\omega_X^{1-p}$ correspond to
homomorphisms
$F_*\O_X\rightarrow \O_{X'}$. A section of $\omega_X^{1-p}$ which
corresponds
to a Frobenius splitting in this way, is called a {\it splitting
section}.
Checking whether a section of $\omega_X^{1-p}$ is a splitting
section can be done locally. More precisely, we have the following
result.

\begin{lem}[\cite{MehtaRamanathan}]
\label{lemmarestriction}
Let $U$ be an open dense subset of a smooth variety $X$.
If a section $s\in \H^0(X, \omega_X^{1-p})$ restricts to a
splitting section $s|_U\in \H^0(U, \omega_U^{1-p})$ on $U$,
then $s$ is a splitting section.
\end{lem}

An immediate consequence of the definition of Frobenius
splitting is

\begin{lem} [\cite{MehtaRamanathan}]
\label{Lemmalinebundle}
Let $X$ be a Frobenius split variety and let $L$ be a
line bundle on $X$ such that $H^i(X,L^m)=0$ for all large $m$
(for a fixed $i$).
Then $H^i(X,L)=0$.
\end{lem}
\begin{proof}
This follows from the fact that if $X$ is Frobenius split
and $L$ is a line bundle on $X$, then there is an injective
map
$$\H^i(X, L) \hookrightarrow \H^i(X, L^p)$$
of abelian groups.
\end{proof}

In particular, Lemma \ref{Lemmalinebundle} implies that
ample line bundles on projective Frobenius split varieties
have vanishing higher cohomology.

\section{Ramification}

In this section we  have the following setup. By
$H=\{ e, \sigma \}$ we  denote the group of order
2, acting nontrivially on a smooth quasiprojective
variety $Y$ over a field of characteristic $p \neq 2$.
We denote the quotient of $Y$ under this action  by $X$
with the corresponding quotient map $\pi$. We will assume,
in addition, that $X$ is smooth.

\begin{lem}
\label{val2}
Assume that $Y=\Spec(B)$ and $X=\Spec(A)$ are affine.
Let $E$ denote an irreducible subvariety of $X$ of
codimension 1 corresponding to a prime ideal $\p$
in $A$. Let $s \in A$ generate $\p$ in the local
ring $A_{\p}$. If $\pi$ is bijective over $E$, then
there exist a unique prime ideal $\p'$ in $B$ over
$\p$. Furthermore, if $v$ denotes the valuation on the
discrete valuation ring $B_{\p'}$, then $v(s)=2$.
\end{lem}
\begin{proof}
Assume that $\p'$ and $\p''$ are two different prime
ideals in $B$ lying over the prime ideal $\p$ in $A$.
Let $E'$ and $E''$ denote the corresponding subvarieties
of $Y$. Then $\sigma(E') = E''$. Choose $y \in E'
\setminus E''$. Since  $\sigma(y) \in E''$, 
$\sigma(y) \neq y$. But $y$ and $\sigma(y)$ both map
to the same point in $E$, which is a contradiction.
This proves the first part of the statement.

Let $E'$  denote the irreducible subvariety
of $Y$ corresponding to the prime ideal $\p'$ in $B$
lying over $\p$. Let $t \in B$ be an element generating
the maximal ideal in the local ring $B_{\p'}$.
Choose $b, b' \in B \setminus \p'$ such that
$$s = t^{v(s)} \frac{b}{b'}.$$
As the product $\sigma(t) t$ is $H$-invariant, we can
find $a,a' \in A \setminus \p$ and a positive integer
$l$ such that
$$ \sigma(t) t = s^l \frac{a}{a'}.$$
Hence, we get
$$s^2 = \sigma(s) s = (\sigma(t) t)^{v(s)} \frac{\sigma(b) b}{\sigma(b')
b'}
= s^{l v(s)}  (\frac{a}{a'})^{v(s)} \frac{\sigma(b) b}{\sigma(b') b'},$$
from which we obtain $l v(s)=2$ (observe that $ \sigma(b) b$ and $
\sigma(b') b' \in  A \setminus \p$).
Assume, if possible,  that $v(s)=1$. Then replacing
$t$ by $t b \sigma(b')$ and $s$ by $s b' \sigma(b')$, we can assume that
$s=t$. Take a nonzero $f \in B$ such that $\sigma (f) = -f$ (e.g.
$f=g-\sigma(g)$
for an element $g$ not invariant under $H$). Since
$H$ is acting trivially on $E'$, it acts trivially on $B/\p'$ and hence
$f$ belongs to $\p'$ (here we are using the assumption that $p \neq 2$).

Write
$$ f = t^{v(f)} \frac{c}{c'}\,,$$
for $c,c' \in B \setminus \p'$.
Applying  $\sigma$  we get,
$$ \sigma (f) = t^{v(f)} \frac{\sigma (c)}{\sigma(c')}.$$
But, by choice,  $\sigma (f) = -f$ and hence  $\sigma (c \sigma(c')) =
-(c
\sigma(c'))$. In particular, $c \sigma(c') \in \p'$. A contradiction,
proving that
$v(s) = 2.$
\end{proof}

\begin{prop}
\label{2}
Let $E$ be an irreducible reduced divisor of $X$ and assume
that $\pi$ is bijective over $E$. Then there exist a
unique irreducible reduced divisor $E'$ of $Y$ mapping onto
$E$. Furthermore, $\pi^*(\O(E))= \O(2E')$.
\end{prop}
\begin{proof}
That there exists a unique (reduced and irreducible) divisor $E'$ in $Y$
mapping
onto $E$ follows from the corresponding local statement
in Lemma \ref{val2}.
Let $s$ be a section of $\O(E)$ with scheme theoretic divisor of zeros
$(s)_0$ equal to $E$. We want to show that $\pi^*(s)$
has divisor of zeros equal to $2 E'$. But this can
be checked locally, and the local statement  follows
from Lemma \ref{val2}.
\end{proof}

\begin{remark}
The above proposition is false, in general, for $p = 2$  and so is the next 
lemma.
\end{remark}

The following lemma is well known.

\begin{lem}
\label{subsetquotient}
Let $V$ be a closed $H$-invariant subvariety of $Y$.
Then $\pi(V)$ (with the reduced closed subscheme structure) is
the quotient of $V$ by $H$.
(For this lemma, it is  not necessary
to assume $Y$ or $X$ to be smooth.)
\end{lem}

The following is an analogue of Hurwitz theorem.
\begin{prop}
\label{canonical}
Let $E = \{ y \in Y : \sigma(y) = y \} $ denote the fixed
point (reduced) subvariety of the action of $H$ on $Y$. If $E$ is a
(closed)
irreducible divisor in $Y$, then $\pi^*(\omega_X) =
\omega_Y \otimes \O(-E)$.
\end{prop}
\begin{proof}
Let $(d\pi)^n : \pi^*(\omega_X) \rightarrow \omega_Y$ denote the
$n$-th ( where $n :=dim(Y)$) exterior power of the differential
$d\pi : \pi^*(\Omega_X) \rightarrow \Omega_Y$ of $\pi$,
and let $\rho$ denote the corresponding global section
of the line bundle $\omega_Y \otimes \pi^*(\omega_X)^{-1}$.
We want to show that the scheme theoretic divisor of zeros $(\rho)_0$
of $\rho$ is equal to $E$:

Let $U$ denote the complement of $E$ in $Y$. Then $U$
is an open subset of $Y$ on which $H$  acts freely.
The restriction of the quotient map $\pi$ to $U$ is
hence \'etale. In
particular,  the support of $(\rho)_0$ must be contained
in $E$. As $E$ is irreducible and $(\rho)_0$ is effective,
 there exists a non-negative integer $l$
such that $(\rho)_0 = lE$. We have to show that $l=1$:
This can be done locally around a point in $E$, so
we may assume that $X$ and $Y$ are affine with
coordinate rings $A \subset B$ respectively.

By Lemma \ref{subsetquotient},  the image $\pi(E)$ (with
the reduced closed subscheme structure) is isomorphic to $E$.
We may therefore think of $E$ as a closed (irreducible) subvariety of
both $X$ and $Y$ (of codim. $1$). Let $\p$ (resp. $\p'$) denote the
prime ideal of height 1 in $A$ (resp. $B$) corresponding
to $E$. Choose $s \in A$ (resp. $t \in B$) generating
$\p$ (resp. $\p'$) in the local ring $A_{\p}$
(resp. $B_{\p'}$).

By Lemma \ref{val2},  we know that there exist
$b,b' \in B \setminus \p'$ such that
$$s = t^2 \frac{b}{b'}.$$
Replacing  $s$ by $s b' \sigma(b')$ and $b$ by $b \sigma(b')$,
we may assume that $b'=1$. Hence $s = t^2 b$. Now choose a point
$z$ in $E$ such that
\begin{itemize}
\item
$E$ is smooth at $z$.
\item
$b(z) \neq 0$.
\item
$\p$ (resp. $\p'$) is generated by $s$ (resp. $t$) in
the local ring $A_{m_z}$ (resp. $B_{m_z'}$), where
$m_z$ (resp. $m_z'$) is the maximal ideal corresponding
to $z$ in $X$ (resp. $Y$).
\end{itemize}
(Since   all these three conditions are separately valid on dense open
sets
in $E$, such a $z$ indeed exists.) As $E$ (by the choice of $z$)
is smooth at $z$, the local ring $A_{m_z}/{\p} =
B_{m_z'}/{\p'}$ is regular. We can therefore choose
elements $s_2, \dots, s_n \in A$ generating the maximal ideal
in this local ring. Hence $ds \wedge ds_2 \wedge \cdots \wedge
ds_n$ (resp. $dt \wedge ds_2 \wedge \cdots \wedge ds_n$) is a
generator of $\pi^*(\omega_X)$ (resp. $\omega_Y$) at
$z$. Let $c \in B_{m_z'}$ be the element such that
$$db \wedge ds_2 \wedge \cdots \wedge ds_n = c \cdot (dt \wedge
ds_2 \wedge \cdots \wedge ds_n).$$
Then
\begin{equation}
\label{l=3}
ds \wedge ds_2 \wedge \cdots \wedge ds_n = t(ct+2b) \cdot (dt \wedge
ds_2 \wedge \cdots \wedge ds_n).
\end{equation}
Noticing that $ct+2b$ is a unit in $B_{m_z'}$ (by the choice
of $z$), it follows that $l=1$ (since $l$ is the exponent of $t$ on the
right side of Equation (\ref{l=3}) above).
\end{proof}

\begin{remark}
All the results in this section are apparently known, but we did not
find
an appropriate reference. Also one can formulate and prove the analogues
of all these results for $H$ replaced by any finite group $G$, provided
that $p$ is coprime to the order of $G$.
\end{remark}

\section{Crepant resolutions}

In this section $X$ will denote a smooth quasiprojective surface over an
algebraically closed field $k$ of  char.
$p \neq 2$. For any  positive integer $n$, as in \S 1,
let $X^{[n]}$ denote the
Hilbert scheme of $n$ points in $X$ and $\psi :
X^{[n]} \to X^{(n)}$  the Hilbert-Chow morphism.
Whenever $Z$ is a smooth variety, we denote by $\omega_Z$
the canonical bundle on $Z$.

\subsection{A fibre diagram}

As in \S 1, let $\Phi : X^n \rightarrow X^{(n)}$ denote the quotient
map.
 Restrictions of $\Phi $ and $\psi $ to $X^n_* := \Phi^{-1}
(X^{(n)}_*)$ and $X^{[n]}_*:= \psi^{-1}
(X^{(n)}_*) $ respectively, yields
the fibre product diagram:

\begin{equation*}
\label{diagram}
\begin{CD}
\tilde{X}^n_* @>\tilde{\psi}>> X^n_* \\
@V\tilde{\Phi}VV   @VV{\Phi}V \\
X^{[n]}_* @>{\psi}>> X^{(n)}_* \\
\end{CD}
\end{equation*}

It is well known that $\tilde{X}^n_*$ is the blow-up of
$X^n_*$ along the big diagonals
$\Delta_{ij} := \{ (x_1, \dots, x_n) \in X^n_*  : x_i = x_j
\}$ ($i < j$), and
that the map $\tilde{\Phi}$ is the quotient map by the
induced $S_n$-action (cf. \cite{Fogarty2}, Lemma 4.4). 
Let $\tilde{E}_{ij}$ denote the
exceptional (reduced) divisor in $\tilde{X}^n_*$ corresponding
to the diagonal $\Delta_{ij}$, and let $\tilde{E}$ denote
the union of the $\tilde{E}_{ij}$.
Let $X^{[n]}_{**}$ denote the open subset $\psi^{-1}
(X^{(n)}_{**})$ in $X^{[n]}_*$, and let $E$ denote
the complement of $X^{[n]}_{**}$ in $X^{[n]}_*$
with the reduced scheme structure. The variety
$E$ is called the {\it exceptional locus} of
$X^{[n]}_*$. Clearly, $E$ is the image of $\tilde{E}_{ij}$
under $\tilde{\Phi}$ for any $i < j$. In particular, $E$
is an irreducible variety.

\subsection{Factorization of $\tilde{\Phi}$}

As mentioned above, the map $\tilde{\Phi} : \tilde{X}_*^n
\rightarrow X^{[n]}_*$ is the quotient of a certain
$S_n$-action on $\tilde{X}_*^n$. We may divide this
quotient into two parts. Let $A_n$ be the
alternating (normal) subgroup of $S_n$, and let $H$
denote the quotient $S_n/A_n$. Let $\tilde{X}^{[n]}_*$
denote the quotient of $\tilde{X}^n_*$ by $A_n$,
and let $\tilde{\Phi}_1$ denote the corresponding
quotient map. Clearly $X^{[n]}_*$ is then the quotient
of $\tilde{X}^{[n]}_*$ by $H$, and we denote the
corresponding quotient map by $\tilde{\Phi}_2$.
Then $\tilde{\Phi} = \tilde{\Phi}_2 \circ
\tilde{\Phi}_1$.

\subsubsection{Description of $\tilde{\Phi}_1$ and
$\tilde{\Phi}_2$}

It is easily seen that $A_n$ is acting freely on
$X_*^n$ and hence also on $\tilde{X}_*^n$.
As $\tilde{X}_*^n$ is smooth, this
implies that the quotient $\tilde{X}^{[n]}_*$ is
also smooth, and that the quotient map is \'etale.
In particular, we get

\begin{lem}
\label{Phi1}
$\tilde{\Phi}_1^*( \omega_{\tilde{X}^{[n]}_*} ) =
\omega_{\tilde{X}_*^n}.$
\end{lem}

All the divisors $\tilde{E}_{ij}$ map to the same divisor
$E'$ in $\tilde{X}^{[n]}_*$. Clearly $H$ acts trivially on $E'$,
hence it follows from Lemma \ref{subsetquotient} that (reduced) $E'$
is isomorphic to $E$. We will however keep the notation $E'$
to emphasize that $E'$ is thought of as a subvariety of
$\tilde{X}^{[n]}_*$. By Proposition \ref{2}, we get

\begin{lem}
\label{2E'}
$\tilde{\Phi}_2^*(\O(E)) = \O(2E')$.

\end{lem}

We also need the following similar result.

\begin{lem}
\label{l}
$\tilde{\Phi}_1^*(\O(E')) = \O( \tilde{E})$.

\end{lem}
\begin{proof}
This follows easily since $\tilde{\Phi}_1$ is an  \'etale map (in particular,
a smooth morphism) and the set theoretic inverse image of $E'$ under 
 $\tilde{\Phi}_1$ is exactly
equal to $ \tilde{E}$. 
\end{proof}

Finally, we need the following result which follows
immediately from Proposition \ref{canonical}.

\begin{lem}
\label{-E'}
$\tilde{\Phi}_2^*(\omega_{X^{[n]}_*}) =
\omega_{\tilde{X}^{[n]}_*} \otimes \O(-E')$.
\end{lem}

\subsection{Crepant resolution}

In this section we will prove the following crucial result.

\begin{thm}
\label{crepant}
Let char. $k \neq 2$. Then
 $\psi : {X^{[n]}_{*}} \to X^{(n)}_* $ is
a crepant resolution,  meaning that $X^{(n)}_*$ is Gorenstein
such that its dualizing line bundle $\omega_{X^{(n)}_*}$ 
 pulls back to the canonical bundle $\omega_{X^{[n]}_*}$
on $X^{[n]}_*$ under $\psi$.
\end{thm}

First we need the following  preparatory lemmas.

\medskip

Recall that  two cycles 
$Z = \sum{m_i Z_i}$ and $Y=\sum{n_j Y_j}$ in an irreducible 
scheme $X$ are 
said to meet {\it  properly}  if $codim(Z_i \cap Y_j)= codim(Z_i) 
+ codim(Y_j)$, whenever $m_i$ and $n_j$ are 
non-zero (cf. \cite{Fulton}, \S 11.4).

\begin{lem}
Let $L$ be a line bundle on  any quasiprojective smooth variety $X$ 
defined over an algebraically closed field $k$, and $p_1,p_2, \dots, p_n$ 
be a finite set of points in $X$. Then there exist an open
subset $U$ in $X$ containing $p_1,p_2, \dots,p_n$ such that 
the restriction of $L$ to $U$ is trivial.  
\end{lem}
\begin{proof}
As any line bundle on a smooth variety is the quotient of two
effective line bundles, we may assume that $L$ is effective. 
Let $s$ be a global section of $L$, and let $(s)_0$ denote 
the divisor of zeros of $s$. By the Moving Lemma (\cite{Fulton},
\S 11.4), there exist 
a divisor $Z$ rationally equivalent to $(s)_0$ such 
that $Z$ meets properly with $ \sum{p_i}$. In
other words, the complement $U$ of the support of $Z$ 
contains $p_1, \dots, p_n$. Since  rationally equivalent 
divisors give rise to isomorphic line bundles (cf. \cite{Fulton}, Example 
2.1.1), $L_{\vert U}$ is trivial.  This  proves the lemma.  
\end{proof}

Let now $X$ be a smooth quasiprojective even dimensional 
variety of dimension $m$, and let $\omega$ denote the canonical 
bundle on $X$. Then the canonical bundle on
$X^n$ is isomorphic to $\omega_n :=\otimes_{i=1}^n 
p_i^*(\omega)$, where $p_i$ is the projection $X^n \to X$ on the $i$-th factor.
We regard $\omega_n$  as a $S_n$-equivariant sheaf on $X^n$ in the 
obvious way. The sheaf $\omega_n^{S_n}$ of $S_n$-invariant sections
of $\omega_n$ can then naturally be thought of as a sheaf on 
$X^{(n)}$. We claim

\begin{lem}
The sheaf $\omega_n^{S_n}$ is a line bundle on $X^{(n)}$.
\end{lem}
\begin{proof}
Let $p = (p_1, \dots, p_n)$ be a point of $X^{(n)}$. 
 By the above Lemma, there exist an open
subset $U$ of $X$ containing $p_1, \dots, p_n$ and 
such that the line bundle $\omega_{\vert U}$  is trivial. As the fibre over 
$p$ (under 
the quotient map) is contained in $U^n$ and as the 
assertion of the lemma is local, we may assume 
that $X=U$. In particular,  we can assume that 
$\omega$ is trivial.

Let $dX$ be a generating global section of $\omega$.
Then $dX_n=\boxtimes_{i=1}^n dX$ is a generating global 
section of $\omega_n$. As $dX$ is an even form, the section
$dX_n$ is $S_n$-invariant, and hence also a global generating
section of $\omega_n^{S_n}$. This proves that $\omega_n^{S_n}$ 
is a line bundle.  
\end{proof}

\begin{lem}
\label{descent} As above, let $X$ be a smooth quasiprojective 
even dimensional variety. 
 Then, there exists a unique  line bundle
 $ L $ on $X^{(n)}$ which restricts
to the canonical bundle on $X^{(n)}_{**}$.

In particular, if char. $k \neq 2$,
$X^{(n)}_{*}$ is Gorenstein with the dualizing line
bundle
$L_{\vert {X^{(n)}_{*}}}$. (We denote $L_{\vert {X^{(n)}_{*}}}$ 
 by  $\omega_{{X^{(n)}_{*}}}$.) 
\end{lem}
\begin{proof} 
Taking $L = \omega_n^{S_n}$, the existence of line bundle $L$ follows from the 
above lemma. 
Since the map $\Phi$  restricted to $X^{n}_{**}$ is \'etale, the
canonical
bundle $\omega_{{X^{(n)}_{**}}}$
pulls back to the canonical bundle  $\omega_{{X^{n}_{**}}}$. Hence,
$L$ restricts to  the canonical bundle on $X^{(n)}_{**}$.
 The uniqueness of $L$ follows since the
codimension of  $X^{(n)} \setminus
X^{(n)}_{**}$ in the normal variety $X^{(n)}$ is $ m \geq 2$.

Since $A_n$ acts freely on (smooth)  $X^n_*$, the quotient
$\tilde{X}^{(n)}_*$
is smooth (and hence Cohen-Macaulay). Further,  $X^{(n)}_{*} =
\tilde{X}^{(n)}_*/H$ and hence it is Cohen-Macaulay (since $p \neq 2$).
Now, the
assertion that   $X^{(n)}_{*}$ is Gorenstein, follows from
\cite{KumarNarasimhan}, Lemma (2.7).
\end{proof}

\begin{remark}
\label{aramova}
Let $X$ be  a
normal and Gorenstein variety $X$ of even dimension over an algebraically 
closed field of char. $p$. Then  $X^{(n)}$ is 
  Gorenstein (and normal) provided $ p > n$. (This is a result due to 
Aramova \cite{Aramova}.) To prove this, apply the `descent' lemma 
(cf., e.g.,  \cite{MehtaRamadas}) to the canonical bundle $\omega_{X^n}$ of
the $S_n$-variety $X^n$ to get a line bundle $L$ on $X^{(n)}$. Moreover, 
$L_{\vert U^{(n)}_{**}}$ is the canonical bundle (where $U \subset X$ is the 
smooth 
locus), since $\Phi_{\vert U^n}$ is an \'etale map. But the complement of
$U^{(n)}_{**}$ in  $X^{(n)}$ has codim. $\geq 2$ and  $X^{(n)}$ is 
Cohen-Macaulay. Hence, by \cite{KumarNarasimhan}, Lemma (2.7), $L$ is the 
dualizing line bundle of   $X^{(n)}$. This proves that  $X^{(n)}$ is 
Gorenstein. 
\end{remark}
\medskip
{\it From now on, we revert to the assumption that $X$ is a smooth 
quasiprojective surface and $p \neq 2$. }
\medskip

\begin{lem}
\label{t}
Let $\omega_{X^{(n)}_*}$ be the dualizing line bundle on $X^{(n)}_*$
guaranteed by the above lemma. Then there exist an integer $t$ such that

$$\psi^*(\omega_{X^{(n)}_*}) \backsimeq \omega_{X^{[n]}_*} \otimes \O(t
E).$$
\end{lem}
\begin{proof}
As $\psi$ is an isomorphism over $X^{(n)}_{**}$ and  the
restriction of $\omega_{X^{(n)}_*}$ to $X^{(n)}_{**}$ is isomorphic to
the canonical bundle, we see that
$$(\psi^*(\omega_{X^{(n)}_*}))_{|X^{[n]}_{**}}  = \omega_{X^{[n]}_{**}}.$$

As $E$ is irreducible, this clearly implies the result.
\end{proof}

\begin{lem}
\label{blowup}
The canonical bundle on $\tilde{X}^n_*$ is given by
$$ \omega_{\tilde{X}^n_*} = \tilde{\psi}^*
(\omega_{X^n_*}) \otimes \O(\tilde{E}),$$
where $\omega_{X^n_*}$ denotes the canonical bundle on
$X^n_*$.
\end{lem}
\begin{proof} Follows from \cite
{Hartshorne},  Exercise II.8.5.
\end{proof}

Now we can prove Theorem \ref{crepant}.

\begin{proof}(\it{of Theorem \ref{crepant}}\rm)
Choose $t \in \mathbb{ Z}$  with the property
as given in
Lemma \ref{t}. We need to show that $t=0$: 
By Lemmas \ref{t},
 \ref{Phi1} - \ref{-E'},
we know that
\begin{equation}
\begin{split}
\label{l=4}
\tilde{\Phi}^* (\psi^*(\omega_{X^{(n)}_*})) & =  \tilde{\Phi}^*(
\omega_{X^{[n]}_*} \otimes \O(t E)) \\
& =  \tilde{\Phi}_1^*(\omega_{\tilde{X}^{[n]}_*}
\otimes \O((2t-1) E')) \\
&  =  \omega_{\tilde{X}^n_*} \otimes \O((2t-1)
\tilde{E}). \
\end{split}
\end{equation}
We want to compare this with an alternative way of
calculating the left  side of the equation above.
Since $\psi \circ \tilde{\Phi} =
\Phi \circ \tilde{\psi}$,

\begin{equation}
\label{l=5}
\tilde{\Phi}^* (\psi^*(\omega_{X^{(n)}_*})) =
\tilde{\psi}^*
(\Phi^*(\omega_{X^{(n)}_*})).
\end{equation}
As $\Phi$ is \'etale over $X^{(n)}_{**}$, the canonical
bundle on $X^{(n)}_{**}$ pulls back to the canonical
bundle on $X^{n}_{**}$. In particular, $\Phi^*(\omega_{X^{(n)}_*})$
 restricts to the canonical bundle on $X^{n}_{**}$.
But the complement of $X^{n}_{**}$ in $X^n_*$ has
codimension $ 2$, which forces $\Phi^*(\omega_{X^{(n)}_*})$
to be  the canonical bundle on $X^n_*$
(as $X^n_*$ is smooth, in particular, normal).
By Lemma \ref{blowup}, we therefore get
\begin{equation}
\label{l=6}
\tilde{\psi}^*(\Phi^*(\omega_{X^{(n)}_*})) = \omega_{\tilde{X}^n_*}
\otimes \O(-\tilde{E}).
\end{equation}
Combining $(\ref{l=4}) - (\ref{l=6})$, we get  $(2t-1) = -1$ (since
$ \O(\tilde{E})$ is a nontorsion element of Pic $\tilde{X}^n_*$),
which forces $t$ to be equal to zero
 as desired.
\end{proof}

The following result in char. $0$ is due to Beauville. 

\begin{cor}
\label{beauville}
Let char. $k > n$. Then  $X^{(n)}$ is Gorenstein and
 $\psi : {X^{[n]}} \to X^{(n)} $ is
a crepant resolution.
\end{cor}
\begin{proof} The assertion that $X^{(n)}$ is Gorenstein
follows by the same argument as for $X^{(n)}_*$ (cf. the proof of
Lemma \ref{descent}) . Now, since the codim.
of
$X^{[n]}\setminus X^{[n]}_* $ in $X^{[n]}$ is at least two, the
corollary follows from the above theorem.
\end{proof}

\section{Frobenius splitting of Hilbert schemes}

Let $X$ be a quasiprojective smooth surface over
an algebraically closed field $k$ of positive char. $p $.
In this section we will prove that $X^{[n]}$ is Frobenius
split if $X$ is Frobenius split. First we need

\begin{lem}
\label{frobsymm}
Let $Y$ be a quasiprojective Frobenius split variety over  $k$.
Then the $n$-th symmetric product $Y^{(n)}$ of $Y$ is Frobenius
split.
\end{lem}
\begin{proof}
Let $\sigma : F_* \O_Y \rightarrow \O_{Y'}$ be a Frobenius
splitting of $Y$. Then $\sigma^{\boxtimes n} : F_* \O_{Y^n} \rightarrow
\O_{(Y^n)'}$ is a Frobenius splitting of the $n$-fold
product of $Y$. As $\sigma^{\boxtimes n}$ is equivariant with respect to  the
natural actions of the symmetric group $S_n$, it takes $S_n$-invariant
functions on $Y^n$ to $S_n$-invariant functions on
$(Y^n)'$. As $\O_{Y^{(n)}}$ is the subsheaf of $\O_{Y^n}$
consisting of $S_n$-invariant functions,
 $\sigma^{\boxtimes n}$ induces a Frobenius splitting of
$Y^{(n)}$.
\end{proof}

\begin{thm}
\label{main}
Let $X$ be a quasiprojective Frobenius split smooth surface over
an algebraically closed field $k$ of  char. $p > 2$ .
Then, for any $n \geq 1,$  the  Hilbert scheme $X^{[n]} $ of $n$ points
in $X$ is Frobenius split.

\end{thm}
\begin{proof}
By Lemma \ref{frobsymm}, the $n$-th symmetric product $X^{(n)}$
is Frobenius split. In particular, $X^{(n)}_{**}$ is Frobenius split.
Let $\sigma'$ be a splitting section of $\omega_{X^{(n)}_{**}}^{1-p}$
on $X^{(n)}_{**}$.
Thinking of  $\sigma'$ as a section of $\omega_{X^{(n)}_*}^{1-p}$
over $X^{(n)}_{**}$,
as $X^{(n)}_*$ is normal and codim. of  $X^{(n)}_* \setminus
X^{(n)}_{**}$
in  $X^{(n)}_*$ is two, we can extend $\sigma'$ to a global
section $\sigma$ of $\omega_{X^{(n)}_*}^{1-p}$ over  $X^{(n)}_*$ (cf.
Lemma \ref{descent}).
Consider the section
$\tilde{\sigma} = \psi^*(\sigma)$ of $\psi^*(\omega_{X^{(n)}_*}^{1-p}) =
\omega_{X^{[n]}_*}^{1-p}$ over $X^{[n]}_*$ (cf. Theorem \ref{crepant}),
and extend it to a section $\hat{\sigma}$ of $
\omega_{X^{[n]}}^{1-p}$ over $X^{[n]}$. (This is possible since
$X^{[n]}$
is smooth, in particular, normal and the codim. of  $X^{[n]} \setminus
X^{[n]}_* $ in  $X^{[n]}$ is at least two.)

We claim that $\hat{\sigma}$ is  a
splitting section of  $
\omega_{X^{[n]}}^{1-p}$ over $X^{[n]}$. To see this, it is enough
to prove that the restriction $\hat{\sigma}'$ of
$\hat{\sigma}$ to $X^{[n]}_{**}$ is a
 splitting section over $X^{[n]}_{**}$. But $X^{[n]}_{**}$
is isomorphic to $X^{(n)}_{**}$ under $\psi$, and
moreover $\hat{\sigma}'$ corresponds
to $\sigma'$ under this isomorphism. As $\sigma'$,
by definition,   Frobenius splits  $X^{(n)}_{**}$,
the result follows.
\end{proof}

\begin{cor}
\label{amplevanishing}
Let $X$ be a smooth projective Frobenius split surface
over a field of characteristic $p > 2$, and let $L$ \
be an ample line bundle on the Hilbert
scheme $X^{[n]}$. Then $L$ has vanishing higher
cohomology.
\end{cor}

\begin{remark}
(a) One can use Corollary \ref{amplevanishing} and the Semicontinuity Theorem
to get a similar vanishing result in characteristic 0.
\medskip

(b) As mentioned by  V. Mehta,  the known list
of Frobenius split smooth surfaces includes

\begin{enumerate}
\item
Projective examples : toric surfaces (in particular $\mathbb{P}^1
\times \mathbb{P}^1$ and $\mathbb{P}^2$), minimal
rational surfaces, ordinary K3 surfaces and ordinary
abelian surfaces. Furthermore, if $s$ is a splitting
section of a smooth surface $X$ which vanishes to order
$(p-1)$ along a point $x$ on $X$, then the blow-up of
$X$ along $x$ is also Frobenius split.
\item
Affine examples : any smooth affine surface is Frobenius split. (In
fact, any smooth affine variety is Frobenius split.)
\end{enumerate}

It is furthermore known that any projective surface, with
Kodaira dimension $\geq 1$, is not Frobenius split. Also
non-ordinary K3 and abelian surfaces are not Frobenius
split.
\medskip

(c) If the punctual Hilbert scheme $H^{[n]}$ 
(i.e. the fibre of the 
Hilbert-Chow 
morphism $\psi$ at $(x, \cdots, x)$ for some $x \in X$) has 
$H^i(H^{[n]}, \O_{H^{[n]}}) = 0$ for all $ i > 0 $, then (under the 
assumptions of Theorem \ref{main} ) $\psi$ is a rational resolution. 
In particular,
for any smooth quasi projective surface $X$ (not necessarily 
Frobenius split), $X^{(n)}$ would be Cohen-Macaulay.

\end{remark}

\bibliographystyle{amsplain}

\begin{thebibliography}{10}

\bibitem{Aramova}
A.G.~Aramova, \emph{Symmetric products of Gorenstein varieties},
J. Algebra \textbf{146} (1992), 482--496.

\bibitem{Fogarty1}
J.~Fogarty, \emph{Algebraic families on an algebraic surface},
  Amer. J. Math. \textbf{90} (1968), 511--521.

\bibitem{Fogarty2}
J.~Fogarty, \emph{Algebraic families on an algebraic surface, II, The Picard 
scheme of the punctual Hilbert scheme},
  Amer. J. Math. \textbf{95} (1973), 660--687.

\bibitem{Fulton} W.~Fulton, \emph{Intersection Theory}, Springer, 1998.


\bibitem{Hartshorne}
R. ~Hartshorne, \emph{Algebraic Geometry}, Springer-Verlag.


\bibitem{KumarNarasimhan} S.~Kumar and M.S.~Narasimhan, \emph{
Picard group of the moduli spaces of $G$-bundles}, Math. Annalen
\textbf{308}
(1997), 155--173.

\bibitem{MehtaRamadas} V.~Mehta and T.~Ramadas, \emph{Moduli of vector
bundles,
Frobenius splitting, and invariant theory },
Ann. of Math. \textbf{144} (1996), 269--313.

\bibitem{MehtaRamanathan}
V.~Mehta and A. Ramanathan,  \emph{Frobenius splitting and cohomology
vanishing for Schubert varieties}, Ann. of Math. \textbf{122} (1985),
27--40.


\bibitem{Serre}
J.P.~Serre, \emph{Algebraic Groups and Class Fields},
Springer-Verlag (1988).

\end{thebibliography}

\providecommand{\bysame}{\leavevmode\hbox to3em{\hrulefill}\thinspace}

\end{document}